\documentclass[final,1p,times,margin=1in]{elsarticle}  

\usepackage{amssymb}
 \usepackage{amsthm}
\usepackage{amscd}
\usepackage{amsmath}
\usepackage{amsfonts}
\usepackage{amssymb}
\usepackage{color}
\usepackage{graphicx}
\usepackage{url}
\newtheorem{theorem}{Theorem}[section]

\usepackage[ruled]{algorithm2e}
\usepackage{mathrsfs}
\usepackage{titletoc}
\usepackage{float}
\usepackage{amsmath}

\usepackage{tikz}
\usepackage{caption}
\usepackage{subcaption}
\usepackage{booktabs}
\usepackage{geometry}

\journal{******}

\begin{document}

\begin{frontmatter}

\title{Core detection via Ricci curvature flows on weighted graphs}

\author[ruc]{Juan Zhao}
\ead{zhaojuan0509@ruc.edu.cn}

\author[ruc]{Jicheng Ma}
\ead{2019202433@ruc.edu.cn}

\author[ruc]{Yunyan Yang\corref{cor1}}
\ead{yunyanyang@ruc.edu.cn}

\author[bnu]{Liang Zhao}
\ead{liangzhao@bnu.edu.cn}

\cortext[cor1]{Corresponding author}

\address[ruc]{School of Mathematics, Renmin University of China, Beijing, 100872, China}
\address[bnu]{School of Mathematical Sciences, Key Laboratory of Mathematics and Complex Systems of MOE,\\
Beijing Normal University, Beijing, 100875, China}

\begin{abstract}
Graph Ricci curvature is crucial as it geometrically quantifies network structure. It pinpoints bottlenecks via negative curvature, identifies cohesive communities with positive curvature, and highlights robust hubs. This guides network analysis, resilience assessment, flow optimization, and effective algorithm design.

In this paper, we derived upper and lower bounds for the weights along several kinds of discrete Ricci curvature flows. As an application, we utilized discrete Ricci curvature flows to detect the core subgraph of a finite undirected graph. The novelty of this work has two aspects.
Firstly, along the Ricci curvature flow, the bounds for weights determine the minimum number of iterations required to ensure weights remain between two prescribed positive constants. In particular, for any fixed graph, we conclude weights can not overflow and can not be treated as zero, as long as the iteration does not exceed a certain number of times;
Secondly, it demonstrates that our Ricci curvature flow method for identifying core subgraphs outperforms prior approaches, such as
page rank, degree centrality, betweenness centrality and closeness centrality.
The codes for our algorithms are available at https://github.com/12tangze12/core-detection-via-Ricci-flow.
\end{abstract}


\end{frontmatter}

\titlecontents{section}[0mm]
                       {\vspace{.2\baselineskip}}
                       {\thecontentslabel~\hspace{.5em}}
                        {}
                        {\dotfill\contentspage[{\makebox[0pt][r]{\thecontentspage}}]}
\titlecontents{subsection}[3mm]
                       {\vspace{.2\baselineskip}}
                       {\thecontentslabel~\hspace{.5em}}
                        {}
                       {\dotfill\contentspage[{\makebox[0pt][r]{\thecontentspage}}]}

\setcounter{tocdepth}{2}



\numberwithin{equation}{section}
\section{Introduction}

Curvature measures a manifold's deviation from flatness. The Riemann curvature tensor fully describes this intrinsic bending.
Contracting it yields the Ricci curvature tensor, governing how volumes evolve along geodesics; positive Ricci curvature causes convergence.
Ricci curvature flow, introduced by Hamilton \cite{hamilton1982ricci}, is a geometric evolution equation smoothing out curvature, defined by
\begin{equation*}
\frac{\partial g_{ij}}{\partial t}=-2{\rm Ric}_{ij}.
\end{equation*}
Overtime, it redistributes curvature, aiming for more uniform geometries. Its groundbreaking application was Perelman’s proof of the Poincar\'e conjecture by evolving 3-manifolds towards standard shapes \cite{perelman2002entropy}.

The concept of Ricci curvature has long been extended to graphs, say Forman's Ricci curvature \cite{ref17, ref18}, Ollivier's Ricci curvature \cite{Ollivier1, ref25}, and Lin-Lu-Yau's Ricci curvature \cite{Lin1}. They share geometric properties similar to those of the manifold case. A positive curvature of an edge indicates a tight connection between two vertices, while a negative curvature indicates a loose connection between two vertices. In 2009, it was proposed by Ollivier \cite{ref25} that the Ricci curvature flow on a weighted graph reads as
\begin{equation}\label{r-flow}w_e^\prime(t)=-\kappa_e(t)w_e(t),\end{equation}
where $w_e$ and $\kappa_e$ denote the weight and Ollivier's Ricci curvature on an edge $e$. Intuitively speaking, along such a flow,
tightly connected points will become tighter, while loosely connected points will become looser.
The discrete Forman-Ricci flow was first  applied to network analysis by Weber et al. \cite{B, A}, characterizing complex networks and detecting changes in dynamic datasets. The long-time behavior of such flows was theoretically explored in \cite{C}.
In a 2019 paper \cite{Ni1}, Ni et al. established a community detection algorithm based on a discrete Ollivier's Ricci curvature flow, combined with a topological surgery. A short time later, in \cite{Lai1}, Lai et al. achieved similar results using a normalized Ricci curvature flow based on Lin-Lu-Yau's Ricci curvature.
The use of discrete Ricci flow for identifying network cores (also known as backbones) was introduced by Weber et al. in \cite{D}, where a more technical definition of the network core was also proposed. Further  applications for backbone detection and related problems were studied in \cite{E}.
For theoretical problems on the continuous curvature flow, such as existence, uniqueness, global existence and convergence, we
refer the readers to \cite{bai,Bai-Lin,Li-Munch}.
Recently, Ma and Yang \cite{Ma1,Ma2,Ma3} proposed  a modified (normalized) Ricci curvature flow and piecewise-linear Ricci curvature flows; discrete versions of these flows were also applied to community detection problems.

The purpose of this article is twofold. One is to derive  upper and lower bounds estimation of solutions for
various discrete Ricci curvature flows. In particular, such a bound leads to global existence of the discrete Ricci curvature flows; it also gives the minimal
number of iterations when weights may overflow or may be treated as zero. The other is to apply the discrete Ricci curvature flow to core detection problem. According to \cite{2011, 1999, 1983, D}, a subgraph of a graph is said to be a {\it core} subgraph if it is tightly connected, and if it is removed, then
the topology of the entire graph will undergo drastic changes. In general, a core subgraph is not uniquely determined. Using the discrete
Ricci curvature flow with an appropriate surgery, we shall design an algorithm to find core subgraphs. From experimental results,
our algorithm outperforms previous methods, such as
page rank, degree centrality, betweenness centrality and closeness centrality.

The remaining part of the article is organized as follows.
In Section 2, we provide some notations and main results. In Section 3, we prove our main results. In Section 4, we introduce the concept of core subgraph and give an algorithm to find it trough Ricci curvature flows. In Section 5, we apply our algorithm to three real-world networks, demonstrating the effectiveness of our approach in comparison with classical methods and hypergraph algorithms.

\section{Notations and main results}
Let $G=(V,E,\mathbf{w})$ be a finite connected weighted graph, where $V=\{z_1,z_2,\cdots,z_n\}$ denotes the vertex set,
$E=\{e_1,e_2,\cdots,e_m\}$ denotes the edge set, and
$\mathbf{w}=(w_{e_1},w_{e_2},\cdots,w_{e_m})$ is a weight vector on $E$. The distance between two vertices $z_i$ and $z_j$ is defined as $d(z_i,z_j)=\inf_{\gamma}\sum_{e\in\gamma}w_e$,
 where $\gamma$ is taken over all paths connecting $z_i$ and $z_j$.

 Given $\alpha\in [0,1]$ and $x\in V$, an
$\alpha$-lazy one-step random walk reads as
\begin{equation*}\mu_{x}^\alpha(z)=\left\{\begin{array}{lll}
\alpha&{\rm if}& z=x\\[1.2ex]
(1-\alpha)\frac{w_{xz}}{\sum_{u\sim x}w_{xu}}&{\rm if}& z\sim x\\[1.2ex]
0&&{\rm otherwise}.
\end{array}\right.\end{equation*}
A function $\mu:V\rightarrow[0,+\infty)$ is said to be a probability measure if $\sum_{x\in V}\mu(x)=1$.
Let $\mu_1$ and $\mu_2$ be two probability measures. A coupling between $\mu_1$ and $\mu_2$ is defined as
a map $A:V\times V\rightarrow[0,1]$ satisfying for all $u, v\in V$,
$\sum_{x\in V}A(u,x)=\mu_1(u)$, $\sum_{y\in V}A(y,v)=\mu_2(v)$.
The Wasserstein distance between $\mu_1$ and $\mu_2$ reads as
\begin{equation*}
W(\mu_1,\mu_2)=\inf_A\sum_{u,v\in V}A(u,v)d(u,v),
\end{equation*}
where $A$ is taken over all couplings between  $\mu_1$ and $\mu_2$. On each edge $e=xy$, Ollivier's
Ricci curvature  \cite{ref25} is defined by
\begin{equation*}\kappa_e^\alpha=1-\frac{W(\mu_x^\alpha,\mu_y^\alpha)}{\rho_e};\end{equation*}
while Lin-Lu-Yau's Ricci curvature \cite{Lin1} reads as
\begin{equation*}\kappa_e=\lim_{\alpha\rightarrow 1-0}\frac{\kappa_e^\alpha}{1-\alpha}.\end{equation*}

In view of \cite{Ma1,Ma2}, we consider a discrete Ricci curvature flow
\begin{equation}\label{discrete-1}\left\{\begin{array}{lll}w_e^{(j+1)}=w_e^{(j)}-s\kappa_e^{(j)}\rho_e^{(j)}\\[1.2ex]
w_e^{(j)}>0,\,\,
w_e^{(0)}=w_{0,e}
\end{array}\right.\end{equation}
for all $j\in\mathbb{N}$ and all $e\in E$, where $s>0$ is the step size, $t_j=js$ is the time of the $j$th iteration, $w_e^{(j)}=w_e(t_j)>0$ and $\kappa_e^{(j)}=\kappa_e(t_j)$ are the weight and the Ricci curvature on $e$ at $t_j$.\\

Our first result is stated as follows.
\begin{theorem}\label{thm2.1}
Let $G=(V,E)$ be a finite graph, $E=\{e_1,e_2,\cdots,e_m\}$ be the edge set, and $\mathbf{w}_0$ be the initial weight on $E$.
Then we have the following two conclusions:\\
$(i)$ If $\kappa:E\rightarrow\mathbb{R}$ is Ollivier's Ricci curvature, then for any $0<s<1$, the flow
 (\ref{discrete-1}) has a unique global solution $(w_e^{(j)})_{j\in\mathbb{N}}$ for all $e\in E$; moreover, there holds
 \begin{equation}\label{weight-est-1}
 (1-s)^jw_{0,e}\leq w_e^{(j)}\leq (1+ms)^j\sum_{\tau\in E}w_{0,\tau},\quad\forall j\in\mathbb{N},\,\,\forall e\in E.
 \end{equation}
$(ii)$  If $\kappa:E\rightarrow\mathbb{R}$ is Lin-Lu-Yau's Ricci curvature, then for any $0<s<1/2$, the  flow
 (\ref{discrete-1}) has a unique global solution $(w_e^{(j)})_{j\in\mathbb{N}}$ for all $e\in E$; moreover, there holds
 \begin{equation}\label{weight-est-2}
 (1-2s)^jw_{0,e}\leq w_e^{(j)}\leq (1+2ms)^j\sum_{\tau\in E}w_{0,\tau},\quad\forall j\in\mathbb{N},\,\,\forall e\in E.
 \end{equation}
\end{theorem}

Let us briefly comment on this theorem.
Firstly, both estimates (\ref{weight-est-1}) and (\ref{weight-est-2}) imply the global existence of solutions to (\ref{discrete-1}) respectively.
Secondly, for any two real numbers $\epsilon_0$ and $M$ satisfying $0<\epsilon_0<\min_{e\in E}\{w_{0,e}\}$ and $M\geq (1+ms)^j\sum_{\tau\in E}w_{0,\tau}$,
the inequality
\begin{equation*}
w_e^{(j)}<\epsilon_0\quad{\rm or}\quad  w_e^{(j)}> M
\end{equation*}
implies
\begin{equation*}
j>\frac{\log\frac{\epsilon_0}{\min_{e\in E}w_{0,e}}}{\log(1-s)}\quad{\rm or}\quad
j>\frac{\log\frac{M}{\sum_{\tau\in E}w_{0,\tau}}}{\log(1+ms)}.
\end{equation*}
In particular, setting $\epsilon_0=10^{-7}$, $w_{0,e}=1$ for all $e\in E$, $s=0.01$, $m=100$ and $M=10^7$, we conclude that
no edge weight is smaller than $\epsilon_0$ within 1612 iteration steps, and that no edge weight is larger than $M$ within
16 iteration steps.

Also we consider a discrete quasi-normalized Ricci curvature flow, which was applied to community detection problem in \cite{Ma1},
namely
\begin{equation}\label{norm-1}
\left\{
\begin{array}{lll}
w_e^{(j+1)}=w_e^{(j)}+s\left(-\kappa_e^{(j)}+\frac{\sum_{\tau\in E}\kappa_\tau^{(j)}\rho_\tau^{(j)}}{\sum_{\tau\in E}w_\tau^{(j)}}
\right)\rho_e^{(j)}\\[1.2ex]
w_e^{(j)}>0,\,\,
w_e^{(0)}=w_{0,e},\,\,\forall j\in\mathbb{N},\, e\in E.
\end{array}\right.
\end{equation}

Our second result reads as follows:

\begin{theorem}\label{thm2.2}
Under the assumptions in Theorem \ref{thm2.1}, we have the following two conclusions:\\
$(i)$ If $\kappa:E\rightarrow\mathbb{R}$ is Ollivier's Ricci curvature, then for any $0<s<1/(m+1)$, the flow
 (\ref{norm-1}) has a unique global solution $(w_e^{(j)})_{j\in\mathbb{N}}$ for all $e\in E$; moreover, there holds
 \begin{equation}\label{weight-est-3}
 (1-(m+1)s)^jw_{0,e}\leq w_e^{(j)}\leq (1+ms)^j\sum_{\tau\in E}w_{0,\tau},\quad\forall j\in\mathbb{N},\,\,\forall e\in E.
 \end{equation}
$(ii)$ If $\kappa:E\rightarrow\mathbb{R}$ is Lin-Lu-Yau's Ricci curvature, then for any $0<s<1/(2m+2)$, the flow
 (\ref{norm-1}) has a unique global solution $(w_e^{(j)})_{j\in\mathbb{N}}$ for all $e\in E$; moreover, there holds
 \begin{equation*}
 (1-2(m+1)s)^jw_{0,e}\leq w_e^{(j)}\leq (1+2(m+1)s)^j\sum_{\tau\in E}w_{0,\tau},\quad\forall j\in\mathbb{N},\,\,\forall e\in E.
 \end{equation*}
 \end{theorem}

 The significance of this theorem is completely analogous to that of Theorem \ref{thm2.1}.\\

Next, we investigate discrete versions of the previous Ricci curvature flow in \cite{ref25,Lai1}. To this end, for any real number $\theta>1$,
it is convenient to define a $\mathbf{\theta}$-{\it surgery} on a graph $G=(V,E,\mathbf{w})$ as follows:
if there exists $e\in E$ such that ${w_e}/{\rho_e}> \theta$,
then $e$ is removed from $E$; if there is no such an edge $e$, then $E$ remains unchanged. Clearly the discrete form of Ricci curvature flow
is as follows:
\begin{equation}\label{flow-1}\left\{\begin{array}{lll}w_e^{(j+1)}=w_e^{(j)}-s\kappa_e^{(j)} w_e^{(j)}\\[1.2ex]
w_e^{(j)}>0,\,\, w_e^{(0)}=w_{0,e}.\end{array}\right.\end{equation}

For this discrete Ricci curvature flow, we have the following theorem.
\begin{theorem}\label{thm2.3}
Fix some $\theta>1$.
Under the assumptions in Theorem \ref{thm2.1}, we have the following two conclusions:\\
  $(i)$ if $\kappa$ is Ollivier's Ricci curvature, then for any $0<s<1$,
the flow (\ref{flow-1}) has a unique global solution $(w_e^{(j)})_{j\in\mathbb{N}}$ for all $e\in E$;
moreover, after possible $\theta$-surgery at each iteration, there holds
 \begin{equation}\label{weight-est-5}
 (1-s)^{j}w_{0,e}\leq w_e^{(j)}\leq (1-s+m\theta s)^j\sum_{\tau\in E}w_{0,\tau}.
 \end{equation}
 $(ii)$ if $\kappa$ is Lin-Lu-Yau's Ricci curvature, then for any $0<s<1/2$,
the flow (\ref{flow-1}) has a unique solution $(w_e^{(j)})_{j\in\mathbb{N}}$ for all $e\in E$;
moreover, after possible $\theta$-surgery at each iteration, there holds
 \begin{equation*}
 (1-2s)^jw_{0,e}\leq w_e^{(j)}\leq (1+2m\theta s)^j\sum_{\tau\in E}w_{0,\tau}.
 \end{equation*}
 \end{theorem}

In \cite{Lai1}, Lai, Bai and Lin established a discrete normalized Ricci curvature flow
\begin{equation}\label{dis-norm}
\left\{\begin{array}{lll}
w_e^{(j+1)}=w_e^{(j)}+s\left(-\kappa_e^{(j)}+\frac{\sum_{\tau\in E}\kappa_\tau^{(j)}w_\tau^{(j)}}{\sum_{\tau\in E}w_\tau^{(j)}}
\right)w_e^{(j)}\\[1.2ex]w_e^{(j)}>0,\,\,
w_e^{(0)}=w_{0,e}.\end{array}\right.
\end{equation}


\begin{theorem}\label{thm2.4}
  Fix some $\theta>1$. Under the assumptions in Theorem \ref{thm2.1}, we have the following two conclusions:\\
$(i)$ if $\kappa$ is Ollivier's Ricci curvature,  then for any $0<s<1/(m\theta)$,
 the flow (\ref{dis-norm}) has a unique global solution $(w_e^{(j)})_{j\in\mathbb{N}}$; moreover, up to $\theta$-surgeries, there holds
 \begin{equation}\label{weight-est-7}
 (1-m\theta s)^{j}w_{0,e}\leq w_e^{(j)}\leq \sum_{\tau\in E}w_{0,\tau},\quad\forall j\in\mathbb{N}.
 \end{equation}
 $(ii)$ if $\kappa$ is Lin-Lu-Yau's Ricci curvature, then for any $0<s<1/(m\theta+2)$,
 the flow (\ref{dis-norm}) has a unique global solution $(w_e^{(j)})_{j\in\mathbb{N}}$, and up to $\theta$-surgeries, there holds
 \begin{equation*}
 (1-(m\theta+2)s)^jw_{0,e}\leq w_e^{(j)}\leq \sum_{\tau\in E}w_{0,\tau},\quad\forall j\in\mathbb{N}.
 \end{equation*}
\end{theorem}

In \cite{Ni1}, Ni et al. utilized a discrete Ricci curvature flow
\begin{equation}\label{rho}\left\{\begin{array}{lll}w_e^{(j+1)}=\rho_e^{(j)}-s\kappa_e^{(j)} \rho_e^{(j)}\\[1.2ex]
w_e^{(j)}>0,\,\, w_e^{(0)}=w_{0,e},\end{array}\right.\end{equation}
together with surgeries,
to solve the community detection problem.
Recently, Li-M\"unch \cite{Li-Munch} obtained a convergent and surgical solution of a flow similar to but distinct from (\ref{rho}).
It is worth emphasizing that the edge weights vary dynamically during the flow (\ref{rho}), whereas in \cite{Li-Munch}, the weights
remain fixed and are independent of the evolving distance.\\

Similarly, we have the following:

\begin{theorem}\label{thm2.5}
Fix some $\theta>1$. Under the assumptions in Theorem \ref{thm2.1}, we have the following two conclusions:\\
  $(i)$ if $\kappa:E\rightarrow\mathbb{R}$ is Ollivier's Ricci curvature, then for any $0<s<1$,
 the flow (\ref{rho}) has a unique global solution $(w_e^{(j)})_{j\in\mathbb{N}}$; moreover, up to $\theta$-surgeries, there holds
 \begin{equation}\label{weight-est-9}
 (1-s)^j\theta^{-j}w_{0,e}\leq w_e^{(j)}\leq (1+ms)^jw_{0,e},\quad\forall j\in\mathbb{N}.
 \end{equation}
 $(ii)$ if $\kappa:E\rightarrow\mathbb{R}$ is Lin-Lu-Yau's Ricci curvature, then for any $0<s<1/2$,
 the flow (\ref{rho}) has a unique global solution $(w_e^{(j)})_{j\in\mathbb{N}}$; moreover, up to  $\theta$-surgeries, there holds
 \begin{equation*}
  (1-2s)^j\theta^{-j}w_{0,e}\leq w_e^{(j)}\leq (1+2ms)^jw_{0,e}.
 \end{equation*}
 \end{theorem}

\section{Proofs of main theorems}
In this section, we shall provide uniform estimates for solutions of various discrete Ricci curvature flows. In particular,
we prove Theorems \ref{thm2.1}-\ref{thm2.5}. The key to proof is curvature estimation.\\

{\it Proof of Theorem \ref{thm2.1}.} $(i)$ Given $e=xy\in E$, $\alpha\in[0,1)$, $s\in(0,1)$ and $t_j=js$ for some $j\in\mathbb{N}$, we have by the definition of the Wasserstein distance
 at $t_j$,
\begin{equation*}
W^{(j)}(\mu_x^\alpha,\mu_y^\alpha)\leq \sum_{u,v\in V}A(u,v)\rho^{(j)}(u,v)\leq \sum_{u,v\in V}A(u,v)\sum_{\tau\in E} w_\tau^{(j)}
=\sum_{\tau\in E} w_\tau^{(j)},
\end{equation*}
where $A:V\times V\rightarrow[0,1]$ is a coupling between $\mu_x^\alpha$ and $\mu_y^\alpha$ such that $\sum_{v\in V}A(u,v)=\mu_x^\alpha(u)$, $\sum_{u\in V}A(u,v)=\mu_y^\alpha(v)$ and $\sum_{u,v\in V}A(u,v)=1$. This leads to
\begin{equation*}
\kappa_e^{(j)}=1-\frac{W^{(j)}(\mu_x^\alpha,\mu_y^\alpha)}{\rho_e^{(j)}}\geq 1-\frac{\sum_{\tau\in E} w_\tau^{(j)}}{\rho_e^{(j)}},
\end{equation*}
which together with $\kappa_e^{(j)}\leq 1$ gives
\begin{equation}\label{right}-\sum_{\tau\in E} w_\tau^{(j)}\leq \rho_e^{(j)}-\sum_{\tau\in E} w_\tau^{(j)}\leq \kappa_e^{(j)}\rho_e^{(j)}\leq \rho_e^{(j)}\leq w_e^{(j)}.\end{equation}
Applying (\ref{right}) to the equation (\ref{discrete-1}), we obtain
\begin{equation*}
(1-s)w_e^{(j)}\leq w_e^{(j+1)}\leq w_e^{(j)}+s\sum_{\tau\in E} w_\tau^{(j)}.
\end{equation*}
Hence,
\begin{equation*}
(1-s)w_e^{(j)}\leq w_e^{(j+1)}\leq \sum_{\tau\in E}w_\tau^{(j+1)}\leq (1+ms)\sum_{\tau\in E}w_\tau^{(j)}.
\end{equation*}
It follows that
\begin{equation*}
(1-s)^{j+1}w_{0,e}\leq w_e^{(j+1)}\leq (1+ms)^{j+1}\sum_{e\in E}w_{0,e},
\end{equation*}
which complete the proof of part (i) of the theorem.\\

$(ii)$ Let $\kappa:E\rightarrow \mathbb{R}$ be Lin-Lu-Yau's Ricci curvature.
According to \cite{bai}, one has for all $e\in E$,
\begin{equation*}-\frac{2}{\rho_e}\max_{\tau\in E}w_\tau\leq\kappa_e\leq 2.\end{equation*}
Clearly, we have
\begin{equation*}
-2\rho_e^{(j)}\leq -\kappa_e^{(j)}\rho_e^{(j)}\leq 2 \max_{\tau\in E}w_\tau^{(j)}.
\end{equation*}
As a consequence,
\begin{equation*}
(1-2s)w_e^{(j)}\leq w_e^{(j+1)}=w_e^{(j)}-s\kappa_e^{(j)}\rho_e^{(j)}\leq w_e^{(j)}+2s\sum_{\tau\in E}w_\tau^{(j)}.
\end{equation*}
Therefore,
\begin{equation*}
(1-2s)^{j+1}w_{0,e}\leq w_e^{(j+1)}\leq \sum_{\tau\in E}w_\tau^{(j+1)}\leq (1+2ms)^{j+1}\sum_{\tau\in E}w_{0,\tau}.
\end{equation*}
This ends the proof of the theorem. $\hfill\Box$\\

 {\it Proof of Theorem \ref{thm2.2}.} $(i)$ According to the proof of Theorem \ref{thm2.1}, we have
 \begin{equation*}
 1-\frac{\sum_{\tau\in E}w_\tau^{(j)}}{\rho_e^{(j)}}\leq \kappa_e^{(j)}\leq 1.
 \end{equation*}
 It follows that
 \begin{equation*}
 -\rho_e^{(j)}\leq -\kappa_e^{(j)}\rho_e^{(j)}\leq -\rho_e^{(j)}+\sum_{\tau\in E}w_\tau^{(j)}
 \end{equation*}
 and
 \begin{equation*}
 -m\rho_e^{(j)}\leq \frac{\sum_{\tau\in E}\kappa_\tau^{(j)}\rho_\tau^{(j)}}{\sum_{\tau\in E}w_\tau^{(j)}}\rho_e^{(j)}\leq \rho_e^{(j)}.
 \end{equation*}
 Inserting the above two inequalities to (\ref{norm-1}), we obtain
 \begin{equation*}
 (1-(m+1)s)w_e^{(j)}\leq w_e^{(j+1)}\leq w_e^{(j)}+s\sum_{\tau\in E}w_\tau^{(j)}.
 \end{equation*}
 Hence, if $0<s<1/(m+1)$, then $(w_e^{(j)})$ is a unique solution of (\ref{norm-1}), and thus
 \begin{equation*}
 (1-(m+1)s)^{j+1}w_{0,e}\leq w_e^{(j+1)}\leq \sum_{\tau\in E}w_\tau^{(j+1)}\leq (1+ms)\sum_{\tau\in E}w_\tau^{(j)}\leq (1+ms)^{j+1}\sum_{\tau\in E}w_{0,\tau}.
 \end{equation*}
 This implies (\ref{weight-est-3}) immediately.\\

$(ii)$ If $\kappa:V\times V\rightarrow \mathbb{R}$ is Lin-Lu-Yau's Ricci curvature, then we have
\begin{equation*}
-2\max_{\tau\in E}w_\tau^{(j)}\leq \kappa_e^{(j)}\rho_e^{(j)}\leq 2\rho_e^{(j)}
\end{equation*}
and
\begin{equation*}
-2m\rho_e^{(j)}\leq \frac{\sum_{\tau\in E}\kappa_\tau^{(j)}\rho_\tau^{(j)}}{\sum_{\tau\in E}w_\tau^{(j)}}\rho_e^{(j)}\leq 2\rho_e^{(j)}.
\end{equation*}
Hence,
\begin{equation*}
-2(m+1)w_e^{(j)}\leq\left(-\kappa_e^{(j)}+\frac{\sum_{\tau\in E}\kappa_\tau^{(j)}\rho_\tau^{(j)}}{\sum_{\tau\in E}w_\tau^{(j)}}
\right)\rho_e^{(j)}\leq 2w_e^{(j)}+2\sum_{\tau\in E}w_\tau^{(j)}.
\end{equation*}
Inserting this into (\ref{norm-1}), we obtain
\begin{equation*}
(1-2(m+1)s)w_e^{(j)}\leq w_e^{(j+1)}\leq \sum_{\tau\in E}w_\tau^{(j+1)}\leq (1+2(m+1)s)\sum_{\tau\in E}w_\tau^{(j)}.
\end{equation*}
As a consequence, if $0<s<1/(2m+2)$, then $(w_e^{j})$ is a unique solution of (\ref{norm-1}) and
\begin{equation*}
(1-2(m+1)s)^{j+1}w_{0,e}\leq w_e^{(j+1)}\leq (1+2(m+1)s)^{j+1}\sum_{\tau\in E}w_{0,\tau}.
\end{equation*}
This completes the proof of the theorem. $\hfill\Box$\\

{\it Proof of Theorem \ref{thm2.3}.} $(i)$ If $\kappa$ is Ollivier's Ricci curvature, we have
\begin{equation*}
-sw_e^{(j)}\leq -s\kappa_e^{(j)}w_e^{(j)}\leq -sw_e^{(j)}+s\frac{w_e^{(j)}}{\rho_e^{(j)}}\sum_{\tau\in E}w_\tau^{(j)}.
\end{equation*}
Under $\theta$-surgeries, we have
\begin{equation*}
(1-s)w_e^{(j)}\leq w_e^{(j+1)}\leq \sum_{\tau\in E}w_\tau^{(j+1)}\leq (1-s+m\theta s)\sum_{\tau\in E}w_\tau^{(j)}.
\end{equation*} 
As a consequence,
\begin{equation*}
(1-s)^{j}w_{0,e}\leq w_e^{(j+1)}\leq (1-s+m\theta s)^{j}\sum_{\tau\in E}w_{0,\tau}.
\end{equation*} 
This implies (\ref{weight-est-5}) immediately.\\

$(ii)$ If $\kappa$ is Lin-Lu-Yau's Ricci curvature, then we have
\begin{equation*}
-2sw_e^{(j)}\leq -s\kappa_e^{(j)}w_e^{(j)}\leq 2s\frac{w_e^{(j)}}{\rho_e^{(j)}}\max_{\tau\in E}w_\tau^{(j)}.
\end{equation*}
Under $\theta$-surgeries, we have
\begin{equation*}
(1-2s)w_e^{(j)}\leq w_e^{(j+1)}\leq \sum_{\tau\in E}w_\tau^{(j+1)}\leq (1+2m\theta s)\sum_{\tau\in E}w_\tau^{(j)}.
\end{equation*}
As a consequence,
\begin{equation*}
(1-2s)^{j}w_{0,e}\leq w_e^{(j+1)}\leq (1+2m\theta s)^{j}\sum_{\tau\in E}w_{0,\tau}.
\end{equation*}
This completes the proof of the theorem. $\hfill\Box$\\

{\it Proof of Theorem \ref{thm2.4}.} $(i)$ If $\kappa$ is Ollivier's Ricci curvature, then we have the following estimates
\begin{equation*}
-1\leq -\kappa_e^{(j)}\leq -1+\frac{\sum_{\tau\in E}w_\tau^{(j)}}{\rho_e^{(j)}},
\end{equation*}
and
\begin{equation*}
\left(1-\sum_{e\in E}\frac{w_e^{(j)}}{\rho_e^{(j)}}\right)\sum_{\tau\in E}w_\tau^{(j)}\leq \sum_{e\in E}\kappa_e^{(j)}w_e^{(j)}\leq \sum_{e\in E}w_e^{(j)}.
\end{equation*}
Thus,
\begin{equation*}
-\sum_{\tau\in E}\frac{w_\tau^{(j)}}{\rho_\tau^{(j)}}\leq -\kappa_e^{(j)}+\frac{\sum_{\tau\in E}\kappa_\tau^{(j)}w_\tau^{(j)}}{\sum_{\tau\in E}w_\tau^{(j)}}\leq \sum_{\tau\in E}\frac{w_\tau^{(j)}}{\rho_\tau^{(j)}}.
\end{equation*}
Therefore,
\begin{equation*}
\left(1-s\sum_{\tau\in E}\frac{w_\tau^{(j)}}{\rho_\tau^{(j)}}\right)w_e^{(j)}\leq w_e^{(j+1)}\leq \left(1+s\sum_{\tau\in E}\frac{w_\tau^{(j)}}{\rho_\tau^{(j)}}\right)w_e^{(j)}.
\end{equation*}
Under $\theta$-surgeries, we have
\begin{equation*}
(1-m\theta s)^{j+1}w_{0,e}\leq w_e^{(j+1)}\leq (1+m\theta s)^{j+1}w_{0,e}.
\end{equation*}
This leads to (\ref{weight-est-7}), since $\sum_{e\in E}w_e^{(j)}=\sum_{e\in E}w_{0,e}$ for all $j\in\mathbb{N}$.\\

$(ii)$ If $\kappa$ is Lin-Lu-Yau's Ricci curvature, then
\begin{equation*}
\left(1-s\left(2+\sum_{\tau\in E}\frac{w_\tau^{(j)}}{\rho_\tau^{(j)}}\right)\right)w_e^{(j)}\leq w_e^{(j+1)}\leq
\left(1+s\left(2+2\frac{\max_{\tau\in E}w_\tau^{(j)}}{\rho_e^{(j)}}\right)\right)w_e^{(j)}.
\end{equation*}
Under the $\theta$-surgeries, we obtain
\begin{equation*}
(1-(m\theta+2)s)w_e^{(j)}\leq w_e^{(j+1)}\leq\sum_{\tau} w_\tau^{(j+1)}\leq (1+2s+2m\theta s)\sum_{\tau}w_\tau^{(j)}.
\end{equation*}
Hence,
\begin{equation*}
(1-(m\theta+2)s)^{j+1}w_e^{(0)}\leq w_e^{(j+1)}\leq (1+2s+2m\theta s)^{j+1}\sum_{\tau}w_\tau^{(0)}.
\end{equation*}
This completes the proof of the theorem.$\hfill\Box$\\

{\it Proof of Theorem \ref{thm2.5}.} $(i)$ If $\kappa$ is Ollivier's Ricci curvature, then
\begin{equation*}
-1\leq -\kappa_e^{(j)}\leq -1+\frac{\sum_{\tau\in E}w_\tau^{(j)}}{\rho_e^{(j)}}.
\end{equation*}
Under $\theta$-surgeries, we have $\frac{w_e}{\theta} \leq \rho_e$. Hence,
\begin{equation*}
-\rho_e^{(j)}\leq -\kappa_e^{(j)}\rho_e^{(j)}\leq -\rho_e^{(j)}+\sum_{\tau\in E}w_\tau^{(j)}.
\end{equation*}
There holds
\begin{equation*}
\frac{1-s}{\theta}w_e^{(j)}\leq(1-s)\rho_e^{(j)}\leq w_e^{(j+1)}\leq w_e^{j}+s\sum_{\tau\in E} w_\tau^{(j)}.
\end{equation*}
Therefore,
\begin{equation*}
\frac{1-s}{\theta}w_e^{(j)}\leq w_e^{(j+1)}\leq \sum_{\tau\in E}w_\tau^{(j+1)}\leq (1+ms)\sum_{\tau\in E}w_\tau^{(j)}.
\end{equation*}
It follows that
\begin{equation*}
(\frac{1-s}{\theta})^{j+1}w_{0,e}\leq w_e^{(j+1)}\leq (1+ms)^{j+1}\sum_{e\in E}w_{0,e},
\end{equation*}
This implies (\ref{weight-est-9}) immediately.\\

$(ii)$  If $\kappa$ is Lin-Lu-Yau's Ricci curvature, then we have
\begin{equation*}
-2\max_{\tau\in E}w_\tau^{(j)}\leq \kappa_e^{(j)}\rho_e^{(j)}\leq 2\rho_e^{(j)}.
\end{equation*}
Under $\theta$-surgeries, we have $\frac{w_e}{\theta} \leq \rho_e$. Thus,
\begin{equation*}
(1-2s)\rho_e^{(j)}\leq w_e^{(j+1)}\leq w_e^{(j)}+2s\sum_{\tau\in E}w_\tau^{(j)}.
\end{equation*}
There holds
\begin{equation*}
\frac{1-2s}{\theta}w_e^{(j)}\leq(1-2s)\rho_e^{(j)}\leq w_e^{(j+1)}\leq w_e^{j}+2s\sum_{\tau\in E} w_\tau^{(j)}.
\end{equation*}
Therefore,
\begin{equation*}
\frac{1-2s}{\theta}w_e^{(j)}\leq w_e^{(j+1)}\leq \sum_{\tau\in E}w_\tau^{(j+1)}\leq (1+2ms)\sum_{\tau\in E}w_\tau^{(j)}.
\end{equation*}
It follows that
\begin{equation*}
(\frac{1-2s}{\theta})^{j+1}w_{0,e}\leq w_e^{(j+1)}\leq (1+2ms)^{j+1}\sum_{e\in E}w_{0,e},
\end{equation*}
which complete the proof of the theorem.  $\hfill\Box$

\section{Core subgraphs and Ricci curvature flow}
In this section, we clarify the following issues: what constitutes a core subgraph; why Ricci curvature flow can be used to locate core subgraphs; and how Ricci curvature flow is used to find core subgraphs.
Hereafter, since the effect of flows (\ref{norm-1}), (\ref{flow-1}), (\ref{dis-norm}) and (\ref{rho}) is the same, we only concern the flow (\ref{discrete-1}).

\subsection{Core subgraphs and their metrics}

Let $G=(V,E)$ be a connected graph. A graph $G^\prime=(V^\prime,E^\prime)$ is said to be a {\it core subgraph}, if $V^\prime\subset V$,
$E^\prime\subset E$, and $G^\prime$ is tightly connected. Denote the induced subgraph of $V\setminus V^\prime$ by $G^\ast$, and
the number of node pairs $\{u,v\} \subseteq V \setminus V'$ that are still connected in $G^\ast$ by $\xi$. Among others,
two metrics of $G^\prime$ are defined as
\begin{equation*}
r_d=\frac{1}{|V^\prime|}\sum_{x\in V^\prime}\frac{\deg_{G^\prime}(x)}{\deg_G(x)}\end{equation*}
and
\begin{equation*}r_{s}=\frac{1}{\xi} \sum_{\{u,v\}\subset V \setminus V',\, \text{dist}_{G^\ast}(u,v) < \infty} \frac{\text{dist}_{G^\ast}(u,v)}{\text{dist}_G(u,v)},
\end{equation*}
where $\deg_{G^\prime}(x)$ is the number of neighbors of $x$ in the graph $G^\prime$, $\deg_G(x)$ is the number of neighbors of $x$ in the graph $G$, $|V^\prime|$ is the number of nodes in $V^\prime$, $\text{dist}_{G}$ is the graph distance on $G$ (each edge has a length $1$), and $\text{dist}_{G^*}$ is the graph distance on $G^*$. Clearly $0\leq r_d\leq 1$, $r_s\geq 1$. In general, as the metric $r_d$ approaches $1$, the connections in the core subgraph become tighter; the larger the metric $r_s$ is, the more the shortest path of the point pairs in the
residual subgraph $G^\ast$ passes through the core nodes. Let us give an example of core subgraph and calculate its metrics.\\

{\it Example $1$.} {\it $G=(V,E)$, $V=\{x_1,x_2,x_3,x_4,x_5,x_6,x_7\}$, $E=\{x_1x_3, x_3x_7, x_ix_{i+1},1\leq i\leq 6\}$. If we take $G_1=(V_1,E_1)$ as a core subgraph, where
$V_1=\{x_1,x_2,x_3\}$ and $E_1=\{x_1x_2,x_1x_3,x_2x_3\}$, then
$r_d={5}/{6}$, $r_s=\frac{13}{12}$. If we take $G_2=(V_2,E_2)$ as another core subgraph, where $V_2=\{x_3,x_4,x_5,x_6,x_7\}$,
$E_2=\{x_3x_7,x_ix_{i+1}, 3\leq i\leq 6\}$, then $r_d=9/10$, $r_s=1$.}\\

For more details on core subgraphs, we refer the readers to early works \cite{2011, 1999, 1983}.

\subsection{Core detection via Ricci curvature flow}

Noting that Ricci curvature flow makes tightly connected points tighter and loosely connected points looser,
one may find core subgraphs through the flow. We shall explain, why one can use the flow to do this and
how to do it, through several explicit examples as follows.\\

{\it Example $2$.} Fix an initial weighted graph $G=(V,E,\mathbf{w}_0)$ as in Figure \ref{fig-1}, with
the initial weight $w_{0,e}=1$ for any $e\in E$.
Apply one iteration of the flow (\ref{discrete-1}) with step size \(s = 0.1\) and parameter \(\alpha = 0.1\).  Based on the updated weights in the second graph, a topological surgery is performed: the three edges with the highest weights (dashed edges) are removed and
the resulting isolated nodes are eliminated. The final graph retains only the inner triangle, which is chosen as a core subgraph.
It is easy to check $r_d=2/3$, but the metric $r_s$ is invalid.
\begin{figure}[H]
\centering

\begin{minipage}{0.25\textwidth}
\centering
\begin{tikzpicture}[scale=0.9, every node/.style={inner sep=0pt, minimum size=0pt}]
\coordinate (x1) at (90:1);
\coordinate (x2) at (210:1);
\coordinate (x3) at (330:1);
\coordinate (x4) at (90:2);
\coordinate (x5) at (210:2);
\coordinate (x6) at (330:2);

\draw (x1)--(x2) node[pos=0.5, left=4pt] {\small 1};
\draw (x2)--(x3) node[pos=0.5, below=4pt] {\small 1};
\draw (x3)--(x1) node[pos=0.5, right=4pt] {\small 1};
\draw (x1)--(x4) node[pos=0.5, right=4pt] {\small 1};
\draw (x2)--(x5) node[pos=0.5, left=4pt] {\small 1};
\draw (x3)--(x6) node[pos=0.5, right=4pt] {\small 1};

\foreach \point in {x1,x2,x3,x4,x5,x6}{
    \fill (\point) circle (1pt);
}
\node[above left=2pt] at (x1) {\small $x_1$};
\node[left=5pt] at (x2) {\small $x_2$};
\node[right=5pt] at (x3) {\small $x_3$};
\node[above=2pt] at (x4) {\small $x_4$};
\node[left=2pt] at (x5) {\small $x_5$};
\node[right=2pt] at (x6) {\small $x_6$};
\end{tikzpicture}
\end{minipage}%
%
\begin{minipage}{0.05\textwidth}
\centering
\begin{tikzpicture}
\draw[->, thick] (0,0) -- (1,0) node[midway, above] {\scriptsize Ricci flow};
\end{tikzpicture}
\end{minipage}%
%
\begin{minipage}{0.25\textwidth}
\centering
\begin{tikzpicture}[scale=0.9, every node/.style={inner sep=0pt, minimum size=0pt}]
\coordinate (x1) at (90:1);
\coordinate (x2) at (210:1);
\coordinate (x3) at (330:1);
\coordinate (x4) at (90:2);
\coordinate (x5) at (210:2);
\coordinate (x6) at (330:2);

\draw (x1)--(x2) node[pos=0.5, left=4pt] {\small 0.87};
\draw (x2)--(x3) node[pos=0.5, below=4pt] {\small 0.87};
\draw (x3)--(x1) node[pos=0.5, right=4pt] {\small 0.87};
\draw[dashed] (x1)--(x4) node[pos=0.5, right=4pt] {\small 0.90};
\draw[dashed] (x2)--(x5) node[pos=0.5, left=4pt] {\small 0.90};
\draw[dashed] (x3)--(x6) node[pos=0.5, right=4pt] {\small 0.90};

\foreach \point in {x1,x2,x3,x4,x5,x6}{
    \fill (\point) circle (1pt);
}
\node[above left=2pt] at (x1) {\small $x_1$};
\node[left=5pt] at (x2) {\small $x_2$};
\node[right=5pt] at (x3) {\small $x_3$};
\node[above=2pt] at (x4) {\small $x_4$};
\node[left=2pt] at (x5) {\small $x_5$};
\node[right=2pt] at (x6) {\small $x_6$};
\end{tikzpicture}
\end{minipage}%
%
\begin{minipage}{0.05\textwidth}
\centering
\begin{tikzpicture}
\draw[->, thick] (0,0) -- (1,0) node[midway, above] {\scriptsize Surgery};
\end{tikzpicture}
\end{minipage}%
%
\begin{minipage}{0.25\textwidth}
\centering
\begin{tikzpicture}[scale=0.9, every node/.style={inner sep=0pt, minimum size=0pt}]
\coordinate (x1) at (90:1);
\coordinate (x2) at (210:1);
\coordinate (x3) at (330:1);

\draw (x1)--(x2) node[pos=0.5, left=4pt] {\small 0.87};
\draw (x2)--(x3) node[pos=0.5, below=4pt] {\small 0.87};
\draw (x3)--(x1) node[pos=0.5, right=4pt] {\small 0.87};

\foreach \point in {x1,x2,x3}{
    \fill (\point) circle (1pt);
}
\node[above left=2pt] at (x1) {\small $x_1$};
\node[left=5pt] at (x2) {\small $x_2$};
\node[right=5pt] at (x3) {\small $x_3$};
\end{tikzpicture}
\end{minipage}
\caption{Finding a core subgraph}
\label{fig-1}
\end{figure}

{\it Example $3$.} The initial weighted graph $G=(V,E,\mathbf{w}_0)$ is described as in Figure \ref{fig-2}. Apply the flow (\ref{discrete-1}) with $s=0.1$ and $\alpha=0.1$. After 5 iterations, it becomes the second graph. Then the edge with the highest weight $x_3x_4$ is deleted,  and the node set $\{x_j\}_{j=1}^6$ is chosen as a core. Finally, the subgraph induced by $\{x_j\}_{j=1}^6$ is $(V,E)$.
Obviously, $r_d=1$ and $r_s=1$.
\begin{figure}[H]
\centering

\begin{minipage}{0.22\textwidth}
\centering
\begin{tikzpicture}[scale=0.9, every node/.style={inner sep=0pt, minimum size=0pt}]
\coordinate (x1) at (-1, 0.866);
\coordinate (x2) at (-1, -0.866);
\coordinate (x3) at (0, 0);
\coordinate (x4) at (1, 0);
\coordinate (x5) at (2, 0.866);
\coordinate (x6) at (2, -0.866);

\draw (x1) -- (x2) node[pos=0.5, right=4pt] {\small 1};
\draw (x2) -- (x3) node[pos=0.5, right=4pt] {\small 1};
\draw (x3) -- (x1) node[pos=0.5, right=4pt] {\small 1};
\draw (x4) -- (x5) node[pos=0.5, left=4pt] {\small 1};
\draw (x5) -- (x6) node[pos=0.5, left=4pt] {\small 1};
\draw (x6) -- (x4) node[pos=0.5, left=4pt] {\small 1};
\draw (x3) -- (x4) node[pos=0.5, below=4pt] {\small 1};

\foreach \point in {x1,x2,x3,x4,x5,x6}{
    \fill (\point) circle (1pt);
}

\node[above left=2pt] at (x1) {\small $x_1$};
\node[below left=2pt] at (x2) {\small $x_2$};
\node[below=4pt] at (x3) {\small $x_3$};
\node[below=4pt] at (x4) {\small $x_4$};
\node[above right=2pt] at (x5) {\small $x_5$};
\node[below right=2pt] at (x6) {\small $x_6$};
\end{tikzpicture}
\end{minipage}%
\hspace{0.005\textwidth}
\begin{minipage}{0.05\textwidth}
\centering
\begin{tikzpicture}
\draw[->, thick] (0,0) -- (1,0) node[midway, above] {\scriptsize Ricci flow};
\end{tikzpicture}
\end{minipage}%
\begin{minipage}{0.05\textwidth}
\centering
\end{minipage}%
\hspace{0.01\textwidth}
\begin{minipage}{0.22\textwidth}
\centering
\begin{tikzpicture}[scale=0.6, every node/.style={inner sep=0pt, minimum size=0pt}]
\coordinate (x1) at (-1, 0.88);
\coordinate (x2) at (-1, -0.88);
\coordinate (x3) at (0, 0);

\coordinate (x4) at (2.0, 0);
\coordinate (x5) at (3.0, 0.88);
\coordinate (x6) at (3.0, -0.88);

\draw (x1) -- (x2) node[pos=0.5, above left=6pt] {\scriptsize 0.71};
\draw (x2) -- (x3) node[pos=0.4, below=6pt] {\scriptsize 0.81};
\draw (x3) -- (x1) node[pos=0.4, above=6pt] {\scriptsize 0.81};

\draw (x4) -- (x5) node[pos=0.4, above=6pt] {\scriptsize 0.81};
\draw (x5) -- (x6) node[pos=0.5, above right=4pt] {\scriptsize 0.71};
\draw (x6) -- (x4) node[pos=0.4, below=4pt] {\scriptsize 0.81};

\draw[dashed] (x3) -- (x4) node[pos=0.5, below=6pt] {\scriptsize 1.15};

\foreach \point in {x1,x2,x3,x4,x5,x6}{
    \fill (\point) circle (1pt);
}

\node[above left=2pt] at (x1) {\small $x_1$};
\node[below left=2pt] at (x2) {\small $x_2$};
\node[below=4pt] at (x3) {\small $x_3$};
\node[below=4pt] at (x4) {\small $x_4$};
\node[above right=2pt] at (x5) {\small $x_5$};
\node[below right=2pt] at (x6) {\small $x_6$};
\end{tikzpicture}
\end{minipage}%
\begin{minipage}{0.05\textwidth}
\centering
\end{minipage}%
\hspace{0.01\textwidth}
\begin{minipage}{0.05\textwidth}
\centering
\begin{tikzpicture}
\draw[->, thick] (0,0) -- (1,0) node[midway, above] {\scriptsize Surgery};
\end{tikzpicture}
\end{minipage}%
\begin{minipage}{0.22\textwidth}
\centering
\begin{tikzpicture}[scale=0.6, every node/.style={inner sep=0pt, minimum size=0pt}]
\coordinate (x1) at (-1, 0.88);
\coordinate (x2) at (-1, -0.88);
\coordinate (x3) at (0, 0);

\coordinate (x4) at (2.0, 0);
\coordinate (x5) at (3.0, 0.88);
\coordinate (x6) at (3.0, -0.88);
\draw (x1) -- (x2) node[pos=0.5, right=4pt] {\small};
\draw (x2) -- (x3) node[pos=0.5, below right=4pt] {\small };
\draw (x3) -- (x1) node[pos=0.5, above right=4pt] {\small };
\draw (x4) -- (x5) node[pos=0.5, above left=2pt] {\small };
\draw (x5) -- (x6) node[pos=0.5, left=4pt] {\small };
\draw (x6) -- (x4) node[pos=0.5, below left=2pt] {\small };

\foreach \point in {x1,x2,x3,x4,x5,x6}{
    \fill (\point) circle (1pt);
}

\node[above left=2pt] at (x1) {\small $x_1$};
\node[below left=2pt] at (x2) {\small $x_2$};
\node[below=4pt] at (x3) {\small $x_3$};
\node[below=4pt] at (x4) {\small $x_4$};
\node[above right=2pt] at (x5) {\small $x_5$};
\node[below right=2pt] at (x6) {\small $x_6$};
\end{tikzpicture}
\end{minipage}
\caption{Finding a core subgraph}
\label{fig-2}
\end{figure}

{\it Example $4$.} The initial weighted graph $G=(V,E,\mathbf{w}_0)$ is described as in Figure \ref{fig-3}. Apply the flow (\ref{discrete-1}) with $s=0.1$ and $\alpha=0.1$. After one iteration, the edge deletion strategy removes all six edges, and the remaining nodes become isolated (any node with degree less than 2 is removed). Finally, $\{x_0\}$ is chosen as a core subgraph. Clearly $r_d=0$, and
$r_s$ is invalid.

\begin{figure}[H]
\centering
\noindent
\begin{minipage}{0.28\textwidth}
\centering
\begin{tikzpicture}[scale=1.1, every node/.style={inner sep=0pt, minimum size=0pt}]
\coordinate (c) at (0,0);
\fill (c) circle (1pt);
\node[below right=2pt] at (c) {\small $x_0$};
\foreach \i/\name in {0/x1,60/x2,120/x3,180/x4,240/x5,300/x6} {
    \coordinate (\name) at (\i:1);
    \draw (c) -- (\name) node[pos=0.5, above=1.5pt] {\scriptsize 1};
    \fill (\name) circle (1pt);
}
\node at (0:1.2) {\small $x_1$};
\node at (60:1.2) {\small $x_2$};
\node at (120:1.2) {\small $x_3$};
\node at (180:1.2) {\small $x_4$};
\node at (240:1.2) {\small $x_5$};
\node at (300:1.2) {\small $x_6$};
\end{tikzpicture}
\end{minipage}%
\hspace{0.015\textwidth}%
\begin{minipage}{0.06\textwidth}
\centering
\begin{tikzpicture}
\draw[->, thick] (0,0) -- (1,0) node[midway, above] {\scriptsize Ricci flow};
\end{tikzpicture}
\end{minipage}%
\begin{minipage}{0.05\textwidth}
\centering
\end{minipage}%
\hspace{0.015\textwidth}%
\begin{minipage}{0.28\textwidth}
\centering
\begin{tikzpicture}[scale=1.1, every node/.style={inner sep=0pt, minimum size=0pt}]
\coordinate (c) at (0,0);
\fill (c) circle (1pt);
\node[below right=2pt] at (c) {\small $x_0$};
\foreach \angle/\name in {0/x1,60/x2,120/x3,180/x4,240/x5,300/x6} {
    \coordinate (\name) at (\angle:1);
    \fill (\name) circle (1pt);
    \draw[dashed] (c) -- (\name) node[pos=0.5, above=1.55pt] {\scriptsize 0.98};
}
\node at (0:1.2) {\small $x_1$};
\node at (60:1.2) {\small $x_2$};
\node at (120:1.2) {\small $x_3$};
\node at (180:1.2) {\small $x_4$};
\node at (240:1.2) {\small $x_5$};
\node at (300:1.2) {\small $x_6$};
\end{tikzpicture}
\end{minipage}%
\hspace{0.015\textwidth}%
\begin{minipage}{0.06\textwidth}
\centering
\begin{tikzpicture}[scale=1]
\node (A) at (0,0) {};
\node (B) at (1,0) {};
\draw[->, thick] (A) -- (B) node[midway, above] {\scriptsize Surgery};
\end{tikzpicture}
\end{minipage}%
\hspace{0.015\textwidth}%
\begin{minipage}{0.18\textwidth}
\centering
\begin{tikzpicture}[scale=1.1, every node/.style={inner sep=0pt, minimum size=0pt}]
\coordinate (c) at (0,0);
\fill (c) circle (1pt);
\node[below right=1pt] at (c) {\small $\mathbf{x_0}$};
\end{tikzpicture}
\end{minipage}
\caption{Finding a core subgraph}
\label{fig-3}
\end{figure}

{\it Example $5$.} The initial weighted graph $G=(V,E,\mathbf{w}_0)$ is described as in Figure \ref{fig-4}. Apply the flow (\ref{discrete-1}) with $s=0.1$ and $\alpha=0.1$. After one iteration and a surgery, take a core set
$\{x_j\}_{j=1}^8$. This induces the core subgraph is $(V,E)$, and thus $r_d=r_s=1$.

\begin{figure}[H]
\centering

\begin{minipage}{0.25\textwidth}
\centering
\begin{tikzpicture}[scale=1.2]
\coordinate (x0) at (0,0);
\fill (x0) circle (1.2pt);
\node[below right=1pt] at (x0) {\small $x_0$};

\foreach \angle/\a/\b/\ia/\ib/\posA/\posB in {
45/x1/x2/1/2/above right/above left,
135/x3/x4/3/4/above left/left,
225/x5/x6/5/6/below left/below right,
315/x7/x8/7/8/below right/right} {
  \coordinate (\a) at ({cos(\angle)+0.4*sin(\angle)}, {sin(\angle)-0.4*cos(\angle)});
  \coordinate (\b) at ({cos(\angle)-0.4*sin(\angle)}, {sin(\angle)+0.4*cos(\angle)});
  \fill (\a) circle (1pt);
  \fill (\b) circle (1pt);
  \draw (x0) -- (\a) node[midway, xshift=-4pt, font=\scriptsize] {1};
  \draw (x0) -- (\b) node[midway, xshift=4pt, font=\scriptsize] {1};
  \draw (\a) -- (\b) node[midway, below, font=\scriptsize] {1};
  \node[\posA=1pt] at (\a) {\small $x_{\ia}$};
  \node[\posB=1pt] at (\b) {\small $x_{\ib}$};
}
\end{tikzpicture}
\end{minipage}%
%
\begin{minipage}{0.05\textwidth}
\centering
\begin{tikzpicture}
\draw[->, thick] (0,0) -- (1,0) node[midway, above] {\scriptsize Ricci flow};
\end{tikzpicture}
\end{minipage}%
%
\begin{minipage}{0.25\textwidth}
\centering
\begin{tikzpicture}[scale=1.2, every node/.style={font=\small}]
\coordinate (x0) at (0,0);
\fill (x0) circle (1.5pt);
\node[below right=2pt] at (x0) {$x_0$};

\foreach \angle/\a/\b/\ia/\ib/\posA/\posB in {
  45/x1/x2/1/2/above right/above left,
  135/x3/x4/3/4/above left/left,
  225/x5/x6/5/6/below left/below right,
  315/x7/x8/7/8/below right/right}
{
  \coordinate (\a) at ({cos(\angle)*1.1 + 0.3*sin(\angle)}, {sin(\angle)*1.1 - 0.3*cos(\angle)});
  \coordinate (\b) at ({cos(\angle)*1.1 - 0.3*sin(\angle)}, {sin(\angle)*1.1 + 0.3*cos(\angle)});
  \fill (\a) circle (1pt);
  \fill (\b) circle (1pt);
  \draw[dashed, gray] (x0) -- (\a) node[midway, xshift=-5pt, font=\scriptsize] {1.00};
  \draw[dashed, gray] (x0) -- (\b) node[midway, xshift=5pt, font=\scriptsize] {1.00};
  \draw (\a) -- (\b) node[midway, above, font=\scriptsize] {0.94};
  \node[\posA=2pt] at (\a) {$x_{\ia}$};
  \node[\posB=2pt] at (\b) {$x_{\ib}$};
}
\end{tikzpicture}
\end{minipage}%
%
\begin{minipage}{0.05\textwidth}
\centering
\begin{tikzpicture}
\draw[->, thick] (0,0) -- (1,0) node[midway, above] {\scriptsize Surgery};
\end{tikzpicture}
\end{minipage}%
%
\begin{minipage}{0.25\textwidth}
\centering
\begin{tikzpicture}[scale=1.2]
\coordinate (x0) at (0,0);
\fill (x0) circle (1pt);
\node[below right=1pt] at (x0) {\small $\mathbf{x_0}$};

\foreach \angle/\a/\b/\ia/\ib/\posA/\posB in {
45/x1/x2/1/2/above right/above left,
135/x3/x4/3/4/above left/left,
225/x5/x6/5/6/below left/below right,
315/x7/x8/7/8/below right/right} {
  \coordinate (\a) at ({cos(\angle)+0.4*sin(\angle)}, {sin(\angle)-0.4*cos(\angle)});
  \coordinate (\b) at ({cos(\angle)-0.4*sin(\angle)}, {sin(\angle)+0.4*cos(\angle)});
  \fill (\a) circle (1pt);
  \fill (\b) circle (1pt);
  \draw (\a) -- (\b) node[midway, below, font=\scriptsize] {0.94};
  \node[\posA=1pt] at (\a) {\small $x_{\ia}$};
  \node[\posB=1pt] at (\b) {\small $x_{\ib}$};
}
\end{tikzpicture}
\end{minipage}
\caption{Finding a core subgraph}
\label{fig-4}
\end{figure}

\subsection{Algorithms}

The algorithm below details the core detection procedure based on discrete Ricci curvature flow. Initially, edge weights are iteratively updated according to the Ricci curvature flow equation (\ref{discrete-1}) with a fixed step size. After a predetermined number of iterations, edges are sorted by their updated weights in descending order, and a specified proportion of the highest-weight edges is removed. The resulting graph may contain isolated nodes, which are identified and ranked by their original degrees in the initial graph. To preserve the target core size, a subset of isolated nodes with the highest original degrees is reinstated into the set of remaining non-isolated nodes. The union of these nodes forms a preliminary core set. The subgraph induced by this set is then extracted from the original graph, and its largest connected component is selected as the final core subgraph. The pseudocode is presented as follows.\\

\begin{algorithm}[H]
\caption{Core detection via Ricci curvature flow on weighted graphs}
\label{algorithm1}
\KwIn{Weighted graph $G=(V,E,w)$; maximum iteration $N$; edge removal ratio $\tau$; step size $s$; curvature parameter $\alpha$.}
\KwOut{Core subgraph $G^{\prime}$.}

\textbf{Step 1: Ricci curvature flow evolution}\;
\For{$i \gets 0$ \KwTo $N-1$}{
    {
        Update $w_e$ using Ricci curvature flow (\ref{discrete-1}) with $s$ and $\alpha$\;
    }
}

\textbf{Step 2: Edge removal and node classification}\;
Sort edges by $w_e^{(T)}$ in descending order and remove the top $\tau\%$ of edges\;
Let $\mathcal{S}$ denote the set of non-isolated nodes and $\mathcal{I} = V \setminus \mathcal{S}$ be the set of isolated nodes\;

\textbf{Step 3: Candidate core node set selection}\;
Compute original degrees $d_G(v)$ for $v \in \mathcal{I}$ and sort $\mathcal{I}$ accordingly\;
Let $M \gets \lfloor |V| / 2 \rfloor$ and $\mathcal{I}'$ be top $M - |\mathcal{S}|$ nodes in $\mathcal{I}$\;
Define core node set $\mathcal{C} = \mathcal{S} \cup \mathcal{I}'$\;

\textbf{Step 4: Final core extraction}\;
Construct the subgraph $G_{\mathcal{C}} = G[\mathcal{C}]$ induced by $\mathcal{C}$ from the original graph structure\;
Identify all connected components of $G_{\mathcal{C}}$\
and let $G^{\prime}$ be the largest connected component\;

\Return $G^{\prime}$ as the detected core subgraph\;
\end{algorithm}
\vspace{1em}

The computational complexity of the proposed core detection algorithm is primarily dominated by the Ricci curvature flow evolution phase. During each of the $N$ iterations, the algorithm updates the weight of every edge based on its discrete Ricci curvature. 
Recall that Ollivier's Ricci curvature along an edge $(x,y)$ is defined as $\kappa(x,y) = 1 - \frac{W(\mu_x, \mu_y)}{d(x,y)}$, where $W(\mu_x, \mu_y)$ denotes the Wasserstein distance between the local probability measures centered at $x$ and $y$, and $d(x,y)$ is the shortest path distance between  two nodes. The shortest path distance $d(x,y)$ is typically computed using Dijkstra's algorithm, with a time complexity of $\mathcal{O}(|E| \log |V|)$. The computation of the Wasserstein distance for a single edge, which corresponds to solving a small optimal transport problem, has a complexity of $\mathcal{O}(D^3)$, where $D$ denotes the average node degree. Consequently, updating all edge weights in one iteration requires $\mathcal{O}(|E| D^3)$ operations, and the overall complexity of the Ricci curvature flow phase is $\mathcal{O}(N |E| D^3)$.
The subsequent steps include sorting the final edge weights ($\mathcal{O}(|E| \log |E|)$), identifying isolated nodes and computing their degrees ($\mathcal{O}(|V|)$), and sorting up to $|V|$ nodes to construct the final core node set ($\mathcal{O}(|V| \log |V|)$). Extracting the largest connected component from the induced subgraph requires an additional $\mathcal{O}(|V| + |E|)$ time. Therefore, the overall time complexity of the algorithm is $\mathcal{O}(N |E| D^3 + |E| \log |E| + |V| \log |V|)$. In practice, due to the cubic dependence on the average degree $D$, the term $\mathcal{O}(N |E| D^3)$ significantly outweighs the others and determines the computational cost of the algorithm.

For comparison, the baseline methods considered in our experiments have lower time complexities: PageRank computed via power iteration requires $\mathcal{O}(N|E|)$, and betweenness centrality computed using Brandes' algorithm for weighted graphs runs in $\mathcal{O}(|V| \cdot |E| + |V|^2 \log |V|)$. While our approach is computationally heavier due to the $\mathcal{O}(|E| D^3)$ term, it leverages geometric information to achieve superior core quality, as demonstrated in the experimental results.

\section{Experiments}
In this section, we evaluate the effectiveness of our core detection algorithm by comparing it with several baseline methods on three real-world networks.

\subsection{Real-world Datasets}

Basic information for real-world networks are listed in Table~\ref{tab:realworld}.

\begin{table}[H]
\centering
\caption{Structural statistics of the three real-world networks analyzed in this study.}
\label{tab:realworld}
\begin{tabular}{lcccccc}
\toprule
Network & Vertices & Edges  & AvgDeg & Density & Diameter \\
\midrule
Cora   & 2485   & 5069       & 4.08  & 0.002 & 19 \\
Citeseer & 2120  & 3679     & 3.47 & 0.002 & 28 \\
Bio-CE-HT & 2617  & 2985   & 2.28 & 0.001 & 20 \\
\bottomrule
\end{tabular}
\end{table}

The Cora and Citeseer datasets \cite{Yang2016} are widely used citation networks, where nodes represent scientific publications and edges denote citation relationships between them. The Cora dataset consists of 2485 publications and 5069 citation links, while the Citeseer dataset contains 2120 publications and 3679 citation links.
The Bio-CE-HT network \cite{bio} captures gene functional associations in Caenorhabditis elegans. Each node represents a gene, and edges denote predicted functional relationships based on multiple biological data sources. The network comprises 2617 nodes and 2985 edges.

Following Theorem \ref{thm2.1}, which ensures validity for $0<s<1$, we set the step size to $s=0.1$ for the Ricci curvature flow in all experiments.
The edge removal ratio is fixed at $\tau=80\%$, meaning the top 80\% of edges (by final weight) are removed after the final Ricci curvature flow iteration. The number of iterations $N$ and the curvature parameter $\alpha$ are adjusted for each dataset to optimize performance: specifically, $N = 50$ and $\alpha=0.8$ for Cora, $N = 12$ and $\alpha = 0.1$ for Citeseer, and $N = 30$ and $\alpha = 0.8$ for Bio-CE-HT.

\subsection{Sensitivity Analysis}
To evaluate the robustness of the proposed Ricci flow method with respect to parameter choices, we conducted two sets of sensitivity experiments on the Citeseer dataset: one varying the step size \(s\) of the Ricci flow, and another perturbing the initial edge weights.

We fixed all other parameters (iterations \(N=12\), curvature parameter \(\alpha=0.1\), removal ratio \(80\%\)) and varied the step size \(s\) from \(0.05\) to \(0.55\). Table~\ref{tab:step_sensitivity} reports the resulting core size and quality metrics \(r_d\) and \(r_s\). The results show that for step sizes up to \(0.2\), the core size and metrics remain relatively stable; beyond \(0.25\), the core size increases noticeably while \(r_s\) exhibits larger fluctuations. The step size $s=0.1$
used in the experiments lies within the stable region.

\begin{table}[H]
\centering
\begingroup
\caption{Sensitivity analysis of the Ricci flow step size on the Citeseer dataset.}
\label{tab:step_sensitivity}
\begin{tabular}{cccc}
\toprule
Step size & Core size & \( r_d \) & \( r_s \) \\
\midrule
0.05 & 342 & 0.76 & 1.50 \\
0.10 & 343 & 0.75 & 1.67 \\
0.15 & 340 & 0.74 & 1.69 \\
0.20 & 361 & 0.73 & 1.72 \\
0.25 & 408 & 0.76 & 1.85 \\
0.30 & 466 & 0.78 & 1.53 \\
0.35 & 407 & 0.76 & 1.46 \\
0.40 & 473 & 0.78 & 1.53 \\
0.45 & 485 & 0.78 & 1.52 \\
0.50 & 488 & 0.77 & 1.52 \\
0.55 & 549 & 0.79 & 1.24 \\
\bottomrule
\end{tabular}
\endgroup
\end{table}

 To assess the influence of the initial edge weights, we applied multiplicative perturbations with factor \(\epsilon=0.5\) (i.e., each weight multiplied by a random factor uniformly drawn from \([0.5,1.5]\)). We performed five independent perturbation runs and compared the resulting core sets with the baseline (original weights). Table~\ref{tab:sensitivity_detailed} reports the core size, \(r_d\), and \(r_s\) for each run, together with the mean and standard deviation (Std) over the five perturbations. The core quality metrics exhibit very low variability, demonstrating that the method is robust to moderate perturbations of the initial weights. 

\begin{table}[H]
\centering
\begingroup
\caption{Sensitivity analysis of initial edge weights on the Citeseer dataset (5 perturbation runs).}
\label{tab:sensitivity_detailed}
\begin{tabular}{lccc}
\toprule
Run & Core size & \(r_d\) & \(r_s\) \\
\midrule
Baseline (original) & 343 & 0.75 & 1.67 \\
Perturbation 1      & 383 & 0.76 & 1.68 \\
Perturbation 2      & 361 & 0.75 & 1.70 \\
Perturbation 3      & 374 & 0.77 & 1.67 \\
Perturbation 4      & 361 & 0.76 & 1.70 \\
Perturbation 5      & 359 & 0.75 & 1.68 \\
\midrule
Mean (perturbed)    & 368 & 0.76 & 1.69 \\
Std  (perturbed)    & 9.85  & 0.01 & 0.01 \\
\bottomrule
\end{tabular}
\endgroup
\end{table}

These experiments confirm that the proposed approach is stable with respect to both the step size (within a reasonable range) and the initial weight configuration, supporting its practical applicability.

\subsection{Comparison with baseline centrality methods}

To evaluate the effectiveness of Algorithm \ref{algorithm1}, we compare it against four widely used node centrality measures that serve as baselines. These include degree centrality, which measures node importance based on the number of directly connected edges, identifying high-degree nodes as structurally central; betweenness centrality, which quantifies the extent to which a node lies on shortest paths between other nodes, highlighting nodes that act as bridges in the network; closeness centrality, defined as the reciprocal of the average shortest path length from a node to all others, reflecting how efficiently a node can reach the entire network; and page rank, a probabilistic measure that estimates the stationary distribution of a random walk on the graph, assigning higher scores to nodes connected to other highly ranked nodes through an iterative computation. For formal definitions and further conceptual discussion of these centrality measures, we refer the readers to~\cite{freeman, newman2,page}.

For each dataset, we first apply Ricci curvature flow method to extract a core subgraph, recording its size as the target core size.
For each of the four baseline methods, we rank all nodes by their corresponding centrality scores and select a connected group of nodes with the highest scores whose size matches that of the core identified by the Ricci curvature flow method.
For each resulting core subgraph, we compute two structural metrics: the core cohesiveness $r_d$ and the average distance stretch $r_s$ after removing the core subgraph. This experimental design ensures a fair comparison across methods by controlling for the size. The comparison results of core extraction methods on the Cora, Citeseer and Bio-CE-HT datasets are presented in Tables \ref{tab2}-\ref{tab4}, respectively.

\begin{table}[H]
\centering
\caption{Core extraction results on the Cora dataset (\(N=50\), \(\alpha=0.8\), s=0.1): core size, core edges, cohesion \(r_d\), and stretch \(r_s\) for the proposed Ricci flow method versus four centrality-based baselines.}
\label{tab2}
\begin{tabular}{lcccc}
\toprule
Method & \#Core Nodes & \#Core Edges & $r_d$ & $r_s$ \\
\midrule
Ricci Flow              & 894 & 1724 & \textbf{0.80} & \textbf{2.17} \\
Page Rank                & 894 & 1874 & 0.62 & 1.00 \\
Degree Centrality       & 894 & 2066 & 0.68 & 1.02 \\
Betweenness Centrality  & 894 & 1743 & 0.64 & 1.34 \\
Closeness Centrality    & 894 & 1960 & 0.78 & 2.03 \\
\bottomrule
\end{tabular}
\begin{flushleft}
\end{flushleft}
\end{table}

\begin{table}[H]
\centering
\caption{Core extraction results on the Citeseer dataset (\(N=12\), \(\alpha=0.1\), s=0.1): core size, core edges, cohesion \(r_d\), and stretch \(r_s\) for the proposed Ricci flow method versus four centrality-based baselines.}
\label{tab3}
\begin{tabular}{lcccc}
\toprule
Method & \#Core Nodes & \#Core Edges & $r_d$ & $r_s$ \\
\midrule
Ricci Flow              & 343 & 860 & \textbf{0.75} & \textbf{1.67} \\
Page Rank                & 343 & 782 & 0.51 & 1.11 \\
Degree Centrality       & 343 & 925 & 0.58 & 1.17 \\
Betweenness Centrality  & 343 & 531 & 0.54 & 1.56 \\
Closeness Centrality    & 343 & 895 & 0.73 & 1.36 \\
\bottomrule
\end{tabular}
\begin{flushleft}
\end{flushleft}
\end{table}

\begin{table}[H]
\centering
\caption{Core extraction results on the Bio-CE-HT dataset (\(N=30\), \(\alpha=0.8\), s=0.1): core size, core edges, cohesion \(r_d\), and stretch \(r_s\) for the proposed Ricci flow method versus four centrality-based baselines.}
\label{tab4}
\begin{tabular}{lcccc}
\toprule
Method & \#Core Nodes & \#Core Edges & $r_d$ & $r_s$ \\
\midrule
Ricci Flow              & 542 & 621 & \textbf{0.74} & \textbf{1.80} \\
Page Rank                & 543 & 688 & 0.55 & 1.00 \\
Degree Centrality       & 542 & 596 & 0.59 & 1.00 \\
Betweenness Centrality  & 542 & 753 & 0.64 & 1.00 \\
Closeness Centrality    & 542 & 786 & 0.73 & 1.20 \\
\bottomrule
\end{tabular}
\begin{flushleft}
\end{flushleft}
\end{table}

The above results indicate that across all three datasets, the Ricci curvature flow method consistently yields core subgraphs with higher core cohesiveness $r_d$ compared to the four baseline centrality measures. Specifically, on the Cora dataset, Ricci curvature flow achieves the highest cohesiveness value of 0.80, outperforming the closest competitor, closeness centrality, which attains 0.78. Similarly, on the Citeseer and Bio-CE-HT datasets, Ricci curvature flow leads with $r_d$ values of 0.75 and 0.74, respectively, indicating more tightly connected core structures.
In terms of average distance stretch $r_s$, which reflects the increase in shortest path lengths after removing the core, our algorithm also shows competitive performance. The Ricci curvature flow method balances cohesiveness and distance stretch effectively, producing cores that maintain structural integrity while preserving key connectivity characteristics of the original network.

Based on our experiments on the three tested datasets, the Ricci flow-based core detection method tends
to identify core subgraphs more effectively than classical centrality measures. These
datasets are generally sparse and homogeneous, suggesting that the method performs
well under such conditions. However, these observations are empirical, and a comprehensive theoretical characterization of the graph structures that favor Ricci flow-based core detection remains an open problem.

\subsection{Comparison with hypergraph Ricci curvature flow method}

To further evaluate the effectiveness of our core detection method, we compare it with a recent algorithm designed for directed or undirected hypergraphs \cite{Sengupta}. This hypergraph algorithm performs 40 Ricci curvature flow iterations, and after every two iterations, it removes the top 8\% of edges ranked by weight and normalizes the edge weights before proceeding to the next iteration. After convergence, the algorithm extracts up to 2 connected (or weakly connected) components as core sets, aiming to identify structurally meaningful cores inherent to high-order interactions in hypergraphs.
We applied this hypergraph Ricci curvature flow algorithm to our datasets, following the parameter settings used in their experiments. The results are summarized in Tables \ref{tab5}-\ref{tab7}.

\begin{table}[H]
\centering
\caption{Comparison of core extraction results on the Cora dataset between our Ricci flow method and the hypergraph Ricci curvature flow algorithm (40 iterations, 8\% edge removal per two iterations).}
\label{tab5}
\begin{tabular}{lcccc}
\toprule
Method & \#Core Nodes & \#Core Edges & $r_d$ & $r_s$ \\
\midrule
Ricci Flow              & 894 & 1724 & 0.80 & 2.17 \\
Hypergraph Algorithm    & 3   & 3   & 0.89 & 1.00 \\
\bottomrule
\end{tabular}
\end{table}

\begin{table}[H]
\centering
\caption{Comparison of core extraction results on the Citeseer dataset between our Ricci flow method and the hypergraph Ricci curvature flow algorithm (40 iterations, 8\% edge removal per two iterations).}
\label{tab6}
\begin{tabular}{lcccc}
\toprule
Method & \#Core Nodes & \#Core Edges & $r_d$ & $r_s$ \\
\midrule
Ricci Flow               & 343 & 860 & 0.75 & 1.67 \\
Hypergraph Core 1        & 6   & 10  & 0.95 & 1.00 \\
Hypergraph Core 2        & 5   &  7  & 0.90 & 1.00 \\
\bottomrule
\end{tabular}
\end{table}

\begin{table}[H]
\centering
\caption{Comparison of core extraction results on the Bio-CE-HT dataset between our Ricci flow method and the hypergraph Ricci curvature flow algorithm (40 iterations, 8\% edge removal per two iterations).}
\label{tab7}
\begin{tabular}{lcccc}
\toprule
Method & \#Core Nodes & \#Core Edges & $r_d$ & $r_s$ \\
\midrule
Ricci Flow              & 542 & 621 & 0.74 & 1.80 \\
Hypergraph Algorithm    & 2   & 1  & 1.00 & 1.00 \\
\bottomrule
\end{tabular}
\end{table}

While evaluating the hypergraph Ricci curvature flow core detection method, we observe that it consistently extracts extremely small core sets when applied to standard graph datasets. For instance, on the Cora dataset, only 3 core nodes are identified, and similarly, only 5 or 2 nodes are found on the Citeseer and Bio-CE-HT datasets. Such minimal core sizes raise concerns about structural relevance, as highlighted in the hypergraph study~\cite{Sengupta}, which states that a core containing less than 5\% or more than 50\% of nodes is typically not meaningful. Although this method achieves high core density $r_d$, the values are inflated due to the small size of the cores, which form dense microstructures. The stretch ratio $r_s$ remains at 1.00 across all datasets, indicating minimal impact on the residual structure.

In contrast, our Ricci curvature flow method identifies significantly larger and structurally influential cores. On the Cora dataset, 894 core nodes are extracted, achieving $r_d = 0.80$ and $r_s = 2.17$, indicating strong internal connectivity and influence on network structure. Similar trends are observed for Citeseer ($r_d = 0.75$, $r_s = 1.67$) and Bio-CE-HT ($r_d = 0.74$, $r_s = 1.80$). These results demonstrate that our method not only uncovers meaningful core sizes but also balances internal cohesiveness and external structural impact more effectively than the hypergraph method.

\section{Conclusion}
In this paper, we studied discrete Ricci curvature flows on weighted graphs from both theoretical and algorithmic perspectives. On the theoretical side, we established explicit upper and lower bounds for edge weights evolving under various forms of discrete Ricci curvature flows. These bounds ensure that the flows remain well-posed and numerically stable within a finite number of iterations, guaranteeing that weights stay within a reasonable range during the process, a crucial prerequisite for practical applications. On the algorithmic side, we proposed a novel method to identify the core subgraph of a given network by combining a discrete Ricci curvature flow with a topological surgery procedure. Experimental evaluations on real-world networks demonstrate that our method consistently outperforms classical approaches such as PageRank, degree centrality, betweenness centrality, and closeness centrality.

While our theoretical results provide stability guarantees for static graphs, several challenging questions remain open. For instance, extending these theoretical bounds to more complex settings such as time-varying or weighted dynamic graphs would significantly broaden the applicability of the approach. Furthermore, the convergence of discrete Ricci curvature flows to a fixed point or a steady state remains an open problem. Our current bounds only prevent catastrophic divergence; proving convergence would likely require additional structural assumptions on the graph or further constraints on the flow parameters.

\section*{Acknowledgements}
This research is partly supported by the National Natural Science Foundation of China (No. 12271039)
and the Open Project Program (No. K202303) of Key Laboratory of Mathematics and Complex Systems, Beijing Normal University.

\section*{Declarations}

\noindent
\textbf{Data availability}:
All data needed are available freely at https://github.com/12tangze12/core-detection-via-Ricci-flow.

\noindent
\textbf{Conflict of interest}: The authors declared no potential conflicts of interest with respect to the research, authorship, and publication of this article.

\noindent
\textbf{Ethics approval}: The research does not involve humans and/or animals. The authors declare that there are no ethics issues to be approved or disclosed.


\end{document}